\documentclass[12pt]{article}
\usepackage{amsfonts,amsmath,amssymb}
          \def\cB{{\cal B}}          \def\cC{{\cal C}}
                    \def\cF{{\cal F}}
          \def\cH{{\cal H}}          
          \def\cK{{\cal K}}          
          \def\cN{{\cal N}}          
                    \def\cR{{\cal R}}
                    \def\cU{{\cal U}}
                    \def\cX{{\cal X}}
\def\cY{{\cal Y}}          \def\cZ{{\cal Z}}

\def\fB{{\mathfrak B}}
\def\fH{{\mathfrak H}}
\def\fg{{\mathfrak g}}
\def\fgh{{\widehat{\mathfrak g}}}

\newcommand{\CC}{{\mathbb C}}

\newcommand{\ZZ}{{\mathbb Z}}

\newcommand{\elpa}{{{\cal A}_{q,p}}}
\newcommand{\elpb}{{{\cal B}_{q,p,\lambda}}}
\newcommand{\un}{\mbox{1\hspace{-1mm}I}}
\newcommand{\id}{\mbox{id}}
\newcommand{\uq}[1]{{{\cal U}_{#1}}}

\newcommand{\sln}{{\widehat{sl}_{N}}}

\newcommand{\dy}[2]{{{\cal D}Y_{#1}^{#2}}}

\newcommand{\sfrac}[2]{{\textstyle{\frac{#1}{#2}}}}

\newcommand{\hypergeom}[4]{{_{2}F_{1}\left(\begin{array}{c}{#1}\quad{#2}\\
                        {#3}\end{array};{#4}\right)}}

\def\qmbox#1{\qquad\mbox{#1}\quad}
\def\Ad{\mathop{\rm Ad}\nolimits}
\def\diag{\mathop{\rm diag}\nolimits}

\def\QTHA{Quasitriangular Hopf Algebra (QTHA)\def\QTHA{QTHA}}
\def\QTQHA{Quasitriangular Quasi-Hopf Algebra (QTQHA)\def\QTQHA{QTQHA}}

\begin{document}
\pagestyle{empty}

\markright{\dotfill Version du \dotfill \today\dotfill\dotfill }

\begin{center}

{\Large \textsf{Yangian and quantum universal solutions \\[3mm] 
    of 
    Gervais--Neveu--Felder 
    equations
    }}

\vspace{10mm}

{\large D. Arnaudon$^a$, J. Avan$^b$, L. Frappat$^{ac}$, {E}. Ragoucy$^a$}

\vspace{10mm}

\emph{$^a$ Laboratoire d'Annecy-le-Vieux de Physique Th{\'e}orique}

\emph{LAPTH, CNRS, UMR 5108, Universit{\'e} de Savoie}

\emph{B.P. 110, F-74941 Annecy-le-Vieux Cedex, France}

\vspace{7mm}

\emph{$^b$ Laboratoire de Physique Th{\'e}orique et Hautes {\'E}nergies}

\emph{LPTHE, CNRS, UMR 7589, Universit{\'e}s Paris VI/VII}

\emph{4, place Jussieu, B.P. 126, F-75252 Paris Cedex 05, France}

\vspace{7mm}

\emph{$^c$ Member of Institut Universitaire de France}

\end{center}

\vfill
\vfill

\begin{minipage}{17.5cm}
\begin{abstract}
We construct universal Drinfel'd twists defining deformations of Hopf
algebra structures based upon simple Lie algebras and contragredient
simple Lie superalgebras.  
In particular, we obtain deformed and dynamical double Yangians.  
Some explicit realisations as evaluation representations are given for
$sl_N$, $sl(1|2)$ and $osp(1|2)$. 
\end{abstract}
\end{minipage}

\vfill
\vfill

\rightline{LAPTH-843/01}
\rightline{PAR-LPTHE 01-12}
\rightline{math.QA/0104181}
\rightline{April 2001}

\baselineskip=16pt

\newpage

\tableofcontents

\newpage

\pagestyle{plain}
\setcounter{page}{1}

\section{Introduction}
\setcounter{equation}{0}

Several consistent deformations of Yangian algebras have been proposed in
the past years, starting with the scaling limit, defined in \cite{KLP,Ko},
of vertex-type quantum elliptic algebras $A_{q,p}(sl_2)$ \cite{FIJKMY}.
Extension of these scaling limits to face-type (so-called ``dynamical'')
elliptic algebras \cite {Fe,GerNeu}, and clarification of their connections
at the level of evaluation representations, were proposed in \cite{Cla} for
structures based upon the Lie algebra $sl_2$.  Construction of several
deformations of Yangian algebras at the universal level, and understanding
thereof as Drinfel'd twists of the centrally extended double Yangian
$\dy{}{}(sl_{2})$ \cite{KT,Kh}, was achieved in \cite{Y2K}, following the
schemes developed in the elucidation as Drinfel'd twists of face and vertex
affine elliptic algebras based upon $sl_N$ \cite{JKOS} and face finite
quantum algebras based upon any simple (contragredient super) Lie algebra
\cite{ABRR}.  The deformed double Yangians were thus characterised as
\QTQHA.

Our purpose here is first of all to extend these universal constructions to
the case of deformations of the centrally extended double Yangians
$\dy{}{}(\fg)$ where $\fg$ is a simple Lie algebra of type $sl_N$ 
or a contragredient simple Lie superalgebra of type $sl(M|N)$ ($M\neq N$).  
We will also construct, by the same techniques, 
consistent deformations of $\uq{}(\fg)$ and $\uq{q}(\fgh)$, for any
$\fg$.  As in the 
previous case,  the existence of such deformations may be conjectured
from considering 
suitable limits of vertex or face-type elliptic quantum affine algebras in
their Lax matrix formulation \cite{Fe,FIJKMY} as a quadratic algebra
(Yang--Baxter equation or RLL formalism \cite{FRT}).  However this limit
procedure yields $R$-matrices which are only interpreted
as \textit{conjectural} evaluation representations of
\textit{hypothetical} universal $R$-matrices for (quasi)-Hopf type
algebraic structures and one must needs establish its actual existence.
This will be
achieved here by identification of these limits as evaluation
representations of universal $R$-matrices for the deformations by
particular Drinfel'd twists, known as 
``shifted-cocycle'' twists \cite{Bab}, of Hopf algebra structures.  
This construction systematically endows these deformations with a
Gervais--Neveu--Felder type \QTQHA\ structure. It is characterised by a
particular form of the universal Yang--Baxter equation, to be explicited
below.

\medskip

We shall first of all construct a deformation $\dy{r}{}(\fg)$, with
$\fg$ (super) unitary, along the
derivation generator $d$, of the centrally extended double Yangian
$\dy{}{}(\fg)$.  
When $\fg$ is taken to be $sl_{N}$, the evaluation
representation of the $R$-matrix for this QTQHA is identified, up 
to a gauge
transformation, with the scaling limit of the $R$-matrix 
for the vertex-type elliptic quantum
affine algebra $\elpa(\sln)$, obtained by sending $q$, $p$ and the spectral
parameter $z$ to $1$ whilst keeping the ratio of their logarithms as finite
parameters \cite{AAFR}. 
We then give as the simplest illustration of the
superalgebra case the evaluation representation of $\dy{r}{}(sl(1|2))$.

\medskip

We will then propose different deformations of Hopf algebra structures,
this time along the Cartan subalgebra of the underlying (finite) Lie
algebra.  For historical reasons, they are called \textit{dynamical}
deformations.

We first recall the previous construction \cite{ABRR} of a twist of the
finite quantum enveloping algebra $\uq{q}(\fg )$ to the dynamical algebra
${\cal B}_{q,\lambda}(\fg )$.  A consistent semi-classical limit of this
procedure then yields a universal twist of shifted-cocycle form
\cite{Bab,BBB}, from the undeformed enveloping algebra $\uq{}(\fg)$ to a
dynamical deformation $\uq{s}(\fg)$ for any simple Lie (contragredient
super) algebra $\fg$, thereafter evaluated for $\fg=sl_N$,
$\fg=osp(1|2)$  and $\fg=sl(1|2)$.

Using now the Hopf algebra inclusion of $\uq{q}(\fg)$ into $\uq{q}(\fgh)$,
the same twist acts on $\uq{q}(\fgh)$ to yield the \QTQHA\
$\uq{q\lambda}(\fgh)$, the $R$-matrix of which may be obtained in the 
$sl_N$ case (under
an evaluated form) as a trigonometric limit $ p \to 0$ of $R$-matrix
for the elliptic
affine face-type algebra $\elpb(\sln)$.

Using finally the Hopf algebra inclusion of $\uq{}(\fg)$ into the extended
double Yangian $\dy{}{}(\fg)$, the previous twist from $\uq{}(\fg)$ to
$\uq{s}(\fg)$ leads from $\dy{}{}(\fg)$ to the dynamical double Yangian
$\dy{s}{}(\fg)$.  The $R$-matrix of this \QTQHA\ may also be obtained in the 
$sl_N$ case
(under an evaluated form) as the scaling limit of $R$-matrix for
the previous algebra
$\uq{q\lambda}(\sln)$.

\section{General setting\label{sect:general}}
\setcounter{equation}{0}

Let $\fg$ be a simple Lie algebra (or a contragredient simple Lie
superalgebra different from $psl(N|N)$) of rank $r_\fg$, with symmetrised
Cartan matrix $A=(a_{ij})$ and inverse $A^{-1}=(d_{ij})$.  In the
superalgebra case, we denote by $[.]$ its $\ZZ_2$ grading.  We denote by
$\cH$ the Cartan subalgebra of $\fg$ with basis $\{ h_{i} \}$ and dual
basis $\{ h^\vee_{i} \}$.  Let $\Pi^+$ be the set of positive roots of
$\fg$ endowed with a normal ordering $<$, i.e. if
$\alpha,\beta,\alpha+\beta\in \Pi^+$ with $\alpha<\beta $, then
$\alpha<\alpha +\beta <\beta $.  Let $\rho$ be the
half-sum of the positive roots (resp. even positive roots)
for a simple Lie algebra (resp. superalgebra).

\noindent
We consider the corresponding quantum universal enveloping (super)
algebra $\cU_q(\fg)$. It is endowed with a Hopf structure.
We use the following coproduct for the generators related to simple
roots
\begin{eqnarray}
  \Delta(e_i) &=& e_i \otimes 1 + q^{h_i} \otimes e_i \;, \label{eq:copro1} \\
  \Delta(f_i) &=& f_i \otimes q^{-h_i} + 1 \otimes f_i \;, \\
  \Delta(h_i) &=& h_i \otimes 1 + 1 \otimes h_i \;.  \label{eq:copro3}
\end{eqnarray}
With this choice, the corresponding universal $R$-matrix was given (up to
$q\leftrightarrow q^{-1}$) in \cite{KT141}.  \\
We regard the quantum affine universal enveloping (super) algebra
$\cU_q(\fgh)$, with universal $R$-matrix, compatible with eqs.
(\ref{eq:copro1})-(\ref{eq:copro3}), given in \cite{KhoTolb}. 
In that case, the Cartan subalgebra is completed with
the derivation and central charge generators $d$ and $c$ respectively.
\medskip

\noindent
We introduce the double Yangian $\dy{}{}(sl_N)$ following \cite{KT,Kh}
with the generators 
\begin{equation}
  e_i^\pm(u) \equiv \pm\sum_{
    \begin{array}{c} 
      {\scriptstyle k\geq 0} \\[-.9ex] {\scriptstyle k<0 } \\
    \end{array}
    }
  e_{i,k} u^{-k-1}\,,
  \qquad
  f_i^\pm(u) \equiv \pm\sum_{
    \begin{array}{c} 
      {\scriptstyle k\geq 0} \\[-.9ex] {\scriptstyle k<0 } \\
    \end{array}
    }
  f_{i,k} u^{-k-1}\,,
  \qquad
  h_i^\pm(u) \equiv 1 \pm\sum_{
    \begin{array}{c} 
      {\scriptstyle k\geq 0} \\[-.9ex] {\scriptstyle k<0 } \\
    \end{array}
    }
  h_{i,k} u^{-k-1}\,.
  \label{eq:modesefh}
\end{equation}
satisfying relations given in \cite{Kh}. The generators related to non-simple 
roots are derived from suitable combinations of generators related to 
simple roots, using Chevalley type relations.
\\
Its universal $R$-matrix, defined in \cite{KT}, obeys the Yang--Baxter
equation. Denoting 
by $\pi_e$ the evaluation representation of $\dy{}{}(sl_N)$, the Lax matrix 
$L=(\pi_e\otimes \un)(\cR)$ realises an FRT-type formalism of
$\dy{}{}(sl_N)$ 
with an $R$-matrix defined by $R=(\pi_e\otimes \pi_e)(\cR)$.

\noindent
Our constructions regarding the double Yangian deformations will relie upon
the following conjecture~\cite{KT}:
The universal $R$-matrix of the centrally extended double
Yangian $\dy{}{}(\fg)$ is given, for any Lie algebra $\fg$, by the
general formula (5.3) in \cite{KT}.
We will also assume that the evaluation of this $R$-matrix corresponds
to the defining $R$-matrix of the Yangian $\cY(sl_N)$. 

\medskip\noindent
We will consider the double Yangian $\dy{}{}(sl(M|N))$. 
In this case, the universal 
$R$-matrix is supposed to be obtained similarly from the double of the
corresponding super-Yangian.  \\
Again, central extensions of $\dy{}{}(\fg)$ will contain
in addition the derivation $d$ and central charge $c$.

\medskip

\noindent
The classes of \QTQHA\ we consider have the particular property to
possess an $R$-matrix that satisfies the
Dynamical Yang--Baxter equation (or Gervais--Neveu--Felder equation):
\begin{equation}
  \cR_{12}(\lambda+h^{(3)}) \, \cR_{13}(\lambda) \,
  \cR_{23}(\lambda+h^{(1)}) = \cR_{23}(\lambda) \,
  \cR_{13}(\lambda+h^{(2)}) \, \cR_{12}(\lambda) \;.
  \label{eq:DYBE0}
\end{equation}
We denote by $\lambda$ a vector with coordinates $(s_{1},\ldots,s_{r_\fg})$
in the basis $\{ h_{i} \}$ in the finite case, or with coordinates
$(s_{1},\ldots,s_{r_\fg},r,s')$ in the basis $\{ h_{i},d,c \}$ otherwise.
Then $\lambda + h^{(1)}$ is the vector with coordinates $s_{i} +
{h^\vee_{i}}^{(1)}$ in the finite case, and
$(\overline{s}_{i},r+c^{(1)},s')$ in other cases.  The coefficients
$\overline{s}_{i}$ are given by $\overline{s}_{i} = 0$ if all $s_{j}$ are
zero and $\overline{s}_{i} = s_{i} + {h^\vee_{i}}^{(1)}$ otherwise.  \\
In the case of superalgebras, the tensor product is graded:
$(a \otimes b)(c \otimes d) = (-1)^{[b][c]} (ac \otimes bd)$.
\\
In representation the spectral parameter (if any) is explicit and the
dynamical Yang--Baxter takes the forms
\begin{equation}
  R_{12}(z,\lambda+h^{(3)}) \, R_{13}(zz',\lambda) \,
  R_{23}(z',\lambda+h^{(1)}) = R_{23}(z',\lambda) \,
  R_{13}(zz',\lambda+h^{(2)}) \, R_{12}(z,\lambda) \;,
  \label{eq:DYBE}
\end{equation}
or
\begin{equation}
  R_{12}(\beta,\lambda+h^{(3)}) \,
  R_{13}(\beta+\beta',\lambda) \,
  R_{23}(\beta',\lambda+h^{(1)}) =
  R_{23}(\beta',\lambda) \,
  R_{13}(\beta+\beta',\lambda+h^{(2)}) \,
  R_{12}(\beta,\lambda) \;,
  \label{eq:DYBE2}
\end{equation}
depending upon the multiplicative or additive nature of the spectral
parameter.
\\
In the case of superalgebras, the $R$-matrix obtained as the evaluation of
the universal $R$-matrix satisfies the graded Yang--Baxter equation
\begin{equation}
        \left(R_{12}\right)_{i_1i_2}^{j_1j_2}
        \left(R_{13}\right)_{j_1i_3}^{k_1j_3}
        \left(R_{23}\right)_{j_2j_3}^{k_2k_3}
        (-1)^{[i_1][i_2]+[i_3][j_1]+[j_2][j_3]} =
        \left(R_{23}\right)_{i_2i_3}^{j_2j_3}
        \left(R_{13}\right)_{i_1j_3}^{j_1k_3}
        \left(R_{12}\right)_{j_1j_2}^{k_1k_2}
        (-1)^{[i_3][i_2]+[i_1][j_3]+[j_2][j_1]}
        \label{eq:YBEgradue}
\end{equation}
Redefining the $R$-matrix as
\begin{equation}
        \left(\widetilde{R}\right)_{i_1i_2}^{j_1j_2} =
        \left(R\right)_{i_1i_2}^{j_1j_2} (-1)^{[i_1][i_2]}
        \label{eq:signe}
\end{equation}
the $R$-matrix $\widetilde{R}$ satisfies now the ordinary (i.e. non-graded)
Yang--Baxter equation (this notation will be used throughout the paper when
dealing with superalgebras).  \\
In the $osp(1|2)$ case (resp. $sl(1|2)$), the basis vectors $v_1$,
$v_2$, $v_3$ of the 
three-dimensional representation have $\ZZ_2$ gradings 1, 0, 0
(resp. 1, 0, 1).

\section{Deformed double Yangian $\dy{r}{}(\fg)$\label{sect:ddy}}
\setcounter{equation}{0}

\subsection{Universal form}

Our aim is to contruct a Drinfel'd twist $\cF$ from $\dy{}{}(\fg)$ to a
deformed double Yangian $\dy{r}{}(\fg)$, which is thereby endowed with a
\QTQHA\ structure.  \\
Let $\omega=\exp(\frac{2i\pi}{r_\fg + 1})$ be a root of unity.  Following
\cite{Y2K,ABRR}, we postulate for $\cF$ the linear identity inspired by
\cite{BR}:
\begin{equation}
  \label{eq:eqlinDYr}
  \cF \equiv \cF(r) = \Ad(\phi^{-1}\otimes \un) (\cF) \cdot \cC \;,
\end{equation}
with
\begin{eqnarray}
  \phi &=& \omega^{h_{0,\rho}} e^{(r+c)d} \;, \\
  \cC &\equiv& e^{\frac12 (c\otimes d + d\otimes c)} \cR \;.
\end{eqnarray}
We will use the following properties:
\begin{itemize}
\item The operator $d$ in the double Yangian
  $\dy{}{}(\fg)$ is defined by
  $[d,e_\alpha(u)]=\frac{d}{du}e_\alpha(u)$ for any root $\alpha $
  (see \cite{Kh}).
\item The operator $d$ satisfies
  $\Delta(d)=d\otimes 1 + 1\otimes d$.
\item The generator $h_{0,\rho}$ of $\dy{}{}(\fg)$ is such that
  \begin{equation}
    h_{0,\rho} e_\alpha (u) = e_\alpha (u)(h_{0,\rho}+(\rho |\alpha))\;,
    \qquad  h_{0,\rho} f_\alpha (u) = f_\alpha (u)(h_{0,\rho}-(\rho |\alpha
)) \;,
    \qquad [h_{0,\rho}, h_\alpha (u)]=0 \;,
  \end{equation}
  and hence $\tau= Ad\left(\omega^{h_{0,\rho}^{(1)}}\right)$
  is idempotent.
\end{itemize}
Equation (\ref{eq:eqlinDYr}) can be solved by
\begin{equation}
  \label{eq:FprodDYr}
  \cF(r) = \prod_k^{\longleftarrow} \cF_k(r) \;,
  \qquad \cF_k(r) = \phi_1^k \cC_{12}^{-1}\phi_1^{-k} \;.
\end{equation}
The solution of the linear equation (\ref{eq:eqlinDYr}) given by the
infinite product (\ref{eq:FprodDYr}) satisfies the following shifted
cocycle relation:
\begin{equation}
  \cF^{(12)}(r) (\Delta\otimes\id)(\cF(r))
  = \cF^{(23)}\left(r+c^{(1)}\right) (\id\otimes\Delta) (\cF(r)) \;.
\end{equation}
The proof follows the same lines as in \cite{Y2K}, inspired by
\cite{JKOS}.

\noindent
The twist $\cF(r)$ defines a \QTQHA\ denoted $\dy{r}{}(\fg)$ with
$R$-matrix $\cR_{\dy{r}{}}(r) = \cF_{21}(r)\cR_{\dy{}{}}\cF_{12}(r)^{-1}$.

\subsection{In representation for $\fg=sl_N$}

A deformed double Yangian $\dy{r}{}(sl_N)$ with $R$-matrix $\overline{R}$
was obtained in \cite{AAFR} (and in \cite{Ko} for $sl_2$) by taking the
scaling limit of the $RLL$ representation of $\elpa(\sln)$.  We now
characterise $\overline{R}$, up to gauge transformation, as an evaluation
representation of the action of the Drinfel'd twist (\ref{eq:FprodDYr}) on
the $R$-matrix of $\dy{}{}(sl_N)$.  In this way, we identify the structure
in \cite{AAFR} as a \QTQHA.

\subsubsection{Gauge transformation of $\dy{r}{}(sl_N)$}

The $R$-matrix of $\dy{r}{}(sl_N)$ is given by
\begin{equation}
  \overline{R}_{a\,,\,b}^{c\,,\,a+b-c}(u,r) = -\frac{1}{N} \; \rho_{DYr}(u) \;
  \frac{\displaystyle \sin\frac{\pi u}{r}\,\sin\frac{\pi}{r}}
  {\displaystyle \sin\frac{\pi(u+1)}{r}} \;
  \overline{S}_{a\,,\,b}^{c\,,\,a+b-c}(u,r)
\end{equation}
where
\begin{equation}
  \overline{S}_{a\,,\,b}^{c\,,\,a+b-c}(u,r) =
  \frac{\displaystyle \sin\frac{\pi}{Nr} (u+1+(b-a) r)}
  {\displaystyle \sin\frac{\pi}{Nr}(u+(b-c) r)
    \,\sin\frac{\pi}{Nr}(1-(a-c) r)}
\end{equation}
The normalisation factor $\rho_{DYr}(u)$ is defined by
\begin{equation}
  \rho_{DYr}(u) = \displaystyle
  \frac{S_{2}(-u \vert r,N) \,
    S_{2}(1+u \vert r,N)}
  {S_{2}(u \vert r,N) \,
    S_{2}(1-u \vert r,N)} \,,
  \label{eqso}
\end{equation}
where $S_{2}(x \vert \omega_{1},\omega_{2})$ is  Barnes' double sine
function of periods $\omega_{1}$ and $\omega_{2}$.

\noindent
We perform the gauge transformation
\begin{equation}
  \label{eq:Rsansbarre}
  R := (V\otimes V) \overline{R} (V\otimes V)^{-1}
\end{equation}
\begin{equation}
  \label{eq:Ssansbarre}
  S := (V\otimes V) \overline{S} (V\otimes V)^{-1}
\end{equation}
with
\begin{equation}
  \label{eq:jauge}
  V_i^j = N^{-1/2} \omega^{(i-1)j} \;.
\end{equation}
A similar transformation was recently exhibited in \cite{YangZhen}.
(\ref{eq:Rsansbarre}) leads to the following expression of the matrix
elements of $S$
\begin{equation}
  \label{eq:SResultat}
  S_{i_1 i_2}^{j_1 j_2}(u) = \delta_{i_1+i_2}^{j_1+j_2}
  \left(
    \delta_{i_2}^{j_1} \; \Omega_{i_2-i_1}\!\!\left(\frac{u}{r}\right)
    + \delta_{i_2}^{j_2} \; \Omega_{i_2-i_1}\!\!\left(\frac{1}{r}\right)
  \right)
\end{equation}
where the function $\Omega_n(x)$ is defined by
\begin{equation}
  \label{eq:Omega}
  \Omega_n(x) := \sum_{k=0}^{N-1} \omega^{nk} \cot\frac{\pi}{N}(x+k)
\end{equation}
that is
\begin{equation}
  \label{eq:OmegaResultat}
  \Omega_n(x) = N\left(
    \frac{e^{i\pi x}}{\sin \pi x} e^{- 2i\pi n x/N} - i \delta_n^0
  \right)
  \qquad \mbox{for} \;\; n\in\{0,...,N-1\}
\end{equation}
Note that in (\ref{eq:OmegaResultat}) the integer $n$ has to be taken in
the interval $[0,\ldots,N-1]$ since the expression (\ref{eq:OmegaResultat})
is not explicitly $N$-periodic in $n$.  \\
In particular, all the non zero entries of $S$ are
\begin{eqnarray}
  S_{aa}^{aa}(u) &=& N \left( \cot \frac{\pi u}{r} + \cot \frac{\pi}{r}\right)
  \label{eq:Sexplicite1} \\
  S_{ab}^{ab}(u) &=& N \frac{e^{i\pi /r}}{\sin \frac{\pi}{r}} e^{- 2i\pi
    (b-a) /Nr}   \qquad \mbox{for} \;\; b-a\in\{1,...,N-1\}
  \label{eq:Sexplicite2} \\
  S_{ab}^{ba}(u) &=& N \frac{e^{i\pi u/r}}{\sin \frac{\pi u}{r}} e^{- 2i\pi
  (b-a) u /Nr}  \qquad \mbox{for} \;\; b-a\in\{1,...,N-1\}
  \label{eq:Sexplicite3}
\end{eqnarray}

\subsubsection{Twist from $\dy{}{}(sl_N)$ to $\dy{r}{}(sl_N)$}

We start with the $R$-matrix of $\dy{}{}(sl_{N})$: its non-vanishing
matrix elements
$R_{i_{1}i_{2}}^{j_{1}j_{2}}$ are given by ($1 \le a,b \le N$)
\begin{equation}
  R_{aa}^{aa} = \rho_{\dy{}{}}(u)
  \qquad\qquad R_{ab}^{ab} = \rho_{\dy{}{}}(u) \frac{u}{u+1}
  \qquad\qquad
  R_{ab}^{ba} = \rho_{\dy{}{}}(u) \frac{1}{u+1}
  \label{eq:rdy}
\end{equation}
\begin{equation}
  \label{eq:rhody}
  \rho_{\dy{}{}}(u) = \frac{\Gamma_{1}(u|N) \;
    \Gamma_{1}(u+N|N)}{\Gamma_{1}(u+1|N) \; \Gamma_{1}(u+N-1|N)}
\end{equation}
We then solve the evaluated linear equation, that is:
\begin{equation}
  \label{eq:eqlin}
  F(u+r) = (\phi\otimes \un)^{-1} F(u) (\phi\otimes \un) R(u+r)
\end{equation}
where
\begin{equation}
  \label{eq:phi}
  \phi = \diag(\omega^{-a})
\end{equation}
and $F$ is now a $N^2\times N^2$ matrix which reduces to $1 \!\times\! 1$
blocks $F_{aa}^{aa}=1$ and to $2\!\times\!2$ blocks:
\begin{equation}
  \label{eq:F2x2}
  \left(
    \begin{array}{cc}
      F_{ab}^{ab}(u) & F_{ab}^{ba}(u) \\
      F_{ba}^{ab}(u) & F_{ba}^{ba}(u)
    \end{array}
  \right) = \left(
    \begin{array}{cc}
      b_{ab}(u) & c_{ab}(u) \\
      c_{ab}'(u) & b_{ab}'(u)
    \end{array}
  \right)\,.
\end{equation}
The linear equation (\ref{eq:eqlin}) is then equivalent to:
\begin{eqnarray}
  \label{eq:eqdiff}
  && \hspace{-20mm}
  \omega^{b-a} (u+2r+1) \, b_{ab}(u+2r) - (u+r+\omega^{b-a}(u+2r)) \,
  b_{ab}(u+r) + (u+r-1) \, b_{ab}(u) = 0 \\
  && \hspace{-20mm}
  (u+2r+1) \, c_{ab}(u+2r) - (u+r+\omega^{a-b}(u+2r)) \, c_{ab}(u+r)
  + \omega^{a-b} (u+r-1) \, c_{ab}(u) = 0
\end{eqnarray}
together with similar equations for $b_{ab}'(u)$ and $c_{ab}'(u)$,
deduced from (\ref{eq:eqdiff}) by the change
$\omega\rightarrow\omega^{-1}$.
The solution is expressed in term of hypergeometric functions
$ \hypergeom{a}{b}{c}{z}$:
\begin{eqnarray}
  \label{eq:soleqdiff}
  b_{ab}(u) &=& \frac{\Gamma\left(\frac{u-1}{r}+1\right)}
  {\Gamma\left(\frac{u}{r}+1\right)} \;
  \hypergeom{-\frac1r}{\frac{u-1}{r}+1}{\frac{u}{r}+1}{\omega^{b-a}}
  \\
  c_{ab}(u) &=&
-\frac{\omega^{b-a}}{r}\;\frac{\Gamma\left(\frac{u-1}{r}+1\right)}
  {\Gamma\left(\frac{u}{r}+2\right)} \;
  \hypergeom{-\frac1r+1}{\frac{u-1}{r}+1}{\frac{u}{r}+2}{\omega^{b-a}}
\end{eqnarray}
Similarly, $b_{ab}'(u)$ and $c_{ab}'(u)$ are obtained from $b_{ab}(u)$ and
$c_{ab}(u)$ by changing $\omega^{b-a}$ into $\omega^{a-b}$.  \\
The twist $F(u)$ given by the collection of $2 \!\times\!  2$ blocks
(\ref{eq:F2x2}) and $1 \!\times\!  1$ blocks $F_{aa}^{aa}=1$ applied to the
$R$-matrix of $\dy{}{}(sl_N)$, Eq.  (\ref{eq:rdy}), i.e.
\begin{equation}
  \label{eq:FRF}
  R^F(u) = F_{21}(-u) R(u) F_{12}^{-1}(u)
\end{equation}
provides the $R$-matrix $R^F(u)$ of the deformed double Yangian
$\dy{r}{}(sl_{N})$, Eq.  (\ref{eq:Rsansbarre}), the non vanishing entries
of which are expressed in terms of
(\ref{eq:Sexplicite1})-(\ref{eq:Sexplicite3}).  The proof follows by a
direct computation using the properties of the hypergeometric functions
$_{2}F_{1}$.

\subsection{In representation for $\fg=sl(1|2)$}

We conjecture the existence of a double Yangian of $sl(1|2)$ endowed with
a Hopf superalgebra structure.  The evaluated $R$-matrix of this double
Yangian is assumed to have the canonical Yang-type simplest rational form:
\begin{equation}
  \widetilde{R}(u) = \frac{u}{u+1} \; \sum_{a, b} 
  (-1)^{[a][b]} E_{aa} \otimes E_{bb} +
  \frac{1}{u+1} \left( \sum_{a\ne b} E_{ab} \otimes E_{ba}
    + \sum_{a} E_{aa} \otimes E_{aa}
  \right)
  \label{eq:dysl12}
\end{equation}
where the gradation is $[1]=[3]=1$, $[2]=0$ (the conventions used for
$sl(1|2)$ are those of \cite{ACF} with the fermionic basis).
\\
A similar evaluation of the twist (\ref{eq:FprodDYr}) now leads to the
following expression for the $R$-matrix of the \QTQHA\
$\dy{r}{}(sl(1|2))$ (with $N=3$):
\begin{eqnarray}
  \widetilde{R}(u) &=& 
   - \; \frac{\sin\frac{\pi (u-1)}{r}}{\sin\frac{\pi(u+1)}{r}} \; 
  ( E_{11} \otimes E_{11} + E_{33} \otimes E_{33} ) + E_{22} \otimes E_{22} 
  \\
  && +\; \frac{\sin\frac{\pi u}{r}}
  {\sin\frac{\pi(u+1)}{r}}  \; 
  \sum_{a< b} (-1)^{[a][b]} 
  \left(
  e^{ i\pi/r + 2i\pi (a-b) /Nr} 
  E_{aa} \otimes E_{bb} 
  + 
  e^{ - i\pi/r -  2i\pi (b-a) /Nr }
    E_{bb} \otimes E_{aa}\right)
  \nonumber  \\ &&
  + \; \frac{
    \sin\frac{\pi }{r}}{\sin\frac{\pi(u+1)}{r}} 
  \sum_{a< b}
  \left(    
  e^{ i\pi u/r + 2i\pi (a-b)u /Nr} 
  E_{ab} \otimes E_{ba} 
  + 
  e^{ - i\pi u/r -  2i\pi (b-a)u /Nr} 
  E_{ba} \otimes E_{ab} 
  \right)
  \nonumber
  \label{eq:dyrsl12}
\end{eqnarray}

\section{Twist from $\uq{q}(\fg )$ to ${\cal B}_{q,\lambda}(\fg )$: a
summary\label{sect:Bql}}
\setcounter{equation}{0}

\subsection{Universal form}

The universal $R$-matrix for $\uq{q}(\fg )$ takes the following form
\begin{equation}
  {\cal R} = \widehat{{\cal R}} \; {\cal K}
  \qmbox{;}
  \widehat{{\cal R}} = \prod_{\gamma\in\Pi^+}
  \widehat{{\cal R}}_{\gamma} \;
  \label{eq:univsln}
\end{equation}
where the product is ordered with respect to $>$, the reversed normal
order on $\Pi^+$. The objects
$\widehat{{\cal R}}_{\gamma}$ and ${\cal K}$ are given by
\begin{equation}
  \widehat{{\cal R}}_{\gamma} = \exp_{q^2} ( -(q-q^{-1}) e_{\gamma}
  \otimes f_{\gamma} )
  \label{eq:rhat}
\end{equation}
and
\begin{equation}
  {\cal K} = q^{- \sum_{ij} d_{ij} h_{i} \otimes h_{j}} \;.
\end{equation}
$e_{\gamma}$, $f_{\gamma}$
are the root generators and $h_{i}$ the Cartan generators in the
Serre-Chevalley basis. In (\ref{eq:rhat}) the $q$-exponential is defined by
\begin{equation}
  \exp_{q}(x) \equiv \sum_{n \in \ZZ^+} \frac{x^n}{(n)_{q}!}
  \qquad
  \mbox{where} \qquad (n)_{q}!  \equiv (1)_{q} (2)_{q} \ldots (n)_{q}
  \;\;
  \mbox{and} \;\; (k)_{q} \equiv \frac{1-q^k}{1-q}
  \label{eq:qexp}
\end{equation}
Expanding the product formula (\ref{eq:univsln}) with respect to a
Poincar\'e--Birkhoff--Witt basis ordered with~$<$, $\widehat\cR$ reads
\begin{equation}
  \label{eq:Rsomme}
  \widehat\cR = \cR \cK^{-1} =
  \un\otimes\un + \sum_{m\in \cZ^*} \sigma_m \; {\bold e}^m
  \otimes {\bold f}^m
\end{equation}
where $\cZ=\mbox{Map}(\Pi^+,\ZZ^+)$, and
$\cZ^*=\cZ\setminus\{(0,\dots,0)\}$.  The term ${\bold e}^m$, (resp.
${\bold f}^m$) denotes an element of the PBW basis of the deformed
enveloping nilpotent subalgebra $\uq{q}(\cN^+)$ (resp.  $\uq{q}(\cN^-)$).

\medskip

\noindent
In \cite{JKOS,ABRR}, it has been shown that there exists a twist $\cF$ from
$\uq{q}(\fg)$ to ${\cal B}_{q,\lambda}(\fg)$ which can be expressed as:
\begin{equation}
  \cF_{[\uq{q}(\fg) \to {\cal B}_{q,\lambda}(\fg)]} = \cK^{-1}
   \widehat{\cF} \cK \;,
   \qquad
   \widehat{\cF} =
  \prod_{k \ge 1}^{\curvearrowleft} \Ad(\phi \otimes \un)^k
  \left( \widehat{\cR}^{-1} \right) \;,
\end{equation}
with
\begin{equation}
  \label{eq:phiBql}
  \phi \equiv q^\cX \equiv
  q^{\sum_{ij} d_{ij} h_i h_j + 2 \sum_i s_i h_i} \qmbox{} s_i\in\CC
\end{equation}
It obeys the cocycle condition
\begin{equation}
  \label{eq:cocycleBql}
  \cF_{12}(w) (\Delta \otimes \un) \cF(w) = \cF_{23}(wq^{{h^\vee}^{(1)}})
  (\un\otimes\Delta)  \cF(w)
\end{equation}
with $w=(w_1,\dots,w_{r_\fg}) = q^s = (q^{s_1},\dots,q^{s_{r_\fg}})\in
\CC^{r_\fg}$,
$w q^{h^\vee}=(w_1 q^{h^\vee_1},\dots,w_{r_\fg} q^{h^\vee_{r_\fg}})$
and $h^\vee_i = \sum_j d_{ij} h_j$. 
\\
This twist satisfies the linear equation
\begin{equation}
  \label{eq:eqlinBql}
  \cF = \Ad(\phi^{-1} \otimes \un) (\cF) \; \cK^{-1} \widehat{\cR} \cK
\end{equation}
the solution of which is uniquely defined \cite{ABRR} once one imposes
\emph{(i)} $\cF\in (\cU_q(\fB^+)\otimes \cU_q(\fB^-))^c$,
\emph{(ii)} its projection on $(\cU_q(\fH)^{\otimes 2})^c$ is
$\un\otimes\un$ where the superscript $c$ denotes a suitable completion.
Under these assumptions, one therefore obtains
\begin{equation}
  \label{eq:FsommeBql}
  \widehat{\cF} = \un\otimes\un +
  \sum_{\{p,r\}\in (\cZ^*)^2} \varphi_{pr}(w) \; {\bold e}^p \otimes
  {\bold f}^r
\end{equation}
where the $\varphi_{pr}(w)$ belong to
$\CC[[s_1,..,s_{r_\fg},s_1^{-1},..,s_{r_\fg}^{-1},\hbar]]\otimes
(\cU_q(\fH)^{\otimes 2})^c$. They are defined recursively  (using
(\ref{eq:eqlinBql})) by
\begin{equation}
  \label{eq:recurvarphiql}
  \left(
    1 - q^{(-2 {h^\vee}^{(1)} + \gamma_p -s| \gamma_p)}
  \right)
  \varphi_{pr}(w) =
  \sum_{    
    \begin{array}{c} 
      {\scriptstyle k+m=p} \\[-.9ex] {\scriptstyle l+m=r } \\[-.9ex]
      {\scriptstyle m\neq 0 }  \\
    \end{array}
    }
  (-1)^{[l][m]} a_p^{km} b_r^{lm} \sigma_m
  q^{(-2 {h^\vee}^{(1)} + \gamma_k -s| \gamma_k)} \;\varphi_{kl}(w) \;.
\end{equation}
In the above equation, $\gamma_p$ is the element of the root
lattice associated to ${\bold e}^p$.
The scalar product $(.|.)$ is given by $(x|y)\equiv
\sum_{i,j}a_{ij}x_i y_j$.
The numbers $a_p^{km}$ and $b_r^{lm}$
are defined by
\begin{equation}
  {\bold e}^k {\bold e}^m \ =\ \sum_{p\in\cZ} a_{p}^{km} {\bold e}^p \
  \mbox{ and }\
  {\bold f}^l {\bold f}^m \ =\ \sum_{r\in\cZ} b_{r}^{lm} {\bold f}^r \;.
\end{equation}
The formulas (\ref{eq:FsommeBql}), (\ref{eq:recurvarphiql}) taken from
\cite{ABRR} will be used in Section \ref{sect:Us} in their $q\rightarrow
1$ limit.

\subsection{In representation for $\fg=sl_N$}

In the fundamental representation for $\fg=sl_N$, we get
\begin{equation}
  R = q^{1/N} \left( \un \otimes \un + (q^{-1}-1) \sum_{a} E_{aa}
    \otimes E_{aa} + (q^{-1}-q) \sum_{a<b} E_{ab} \otimes E_{ba} \right)
  \label{eq:reprsln}
\end{equation}
the $N \!\times\! N$ matrices $E_{ab}$ being the usual elementary matrices
with entry 1 in position $(a,b)$ and 0 elsewhere.
\\
The twist is represented by
\begin{eqnarray}
  \widehat{F} &=& \prod_{k \ge 1}^{\curvearrowleft}
  (B \otimes \un)^k \; \widehat{R}^{-1} \; (B
  \otimes \un)^{-k}
  \\
  F_{[\uq{q}(sl_{N}) \to {\cal B}_{q,\lambda}(sl_{N})]}
  &=& \un \otimes \un + (q-q^{-1}) \sum_{a<b} \frac{w_{ab}}{1-w_{ab}}
  \; E_{ab} \otimes E_{ba}
  \label{eq:freprbqn}
\end{eqnarray}
with
$  B = q^{\frac{N-1}{N}} \; \diag(q^{x_a}) $,
$w_{ab} = q^{x_{a} - x_{b}}$ and
$x_{a} = 2 s_{a} - 2 s_{a-1}$ with $s_{0} = s_{N} = 0$.
\\
The non vanishing elements $R_{i_{1}i_{2}}^{j_{1}j_{2}}$ of the
$R$-matrix of ${\cal B}_{q,\lambda}(sl_{N})$
are then given by ($1 \le a,b \le N$)
\begin{eqnarray}
        && R_{aa}^{aa} =  q^{1/N} \frac{1}{q} \nonumber \\
        && R_{ab}^{ab} = q^{1/N} \left\{
        \begin{array}{ll}
                1 & \mbox{if} \;\; b>a \\
                \displaystyle
                \frac{(1-q^2w_{ab})(1-q^{-2}w_{ab})}{(1-w_{ab})^2} &
                \mbox{if} \;\; b<a \\
        \end{array}
        \right.  \\
        && R_{ab}^{ba} = q^{1/N} (q-q^{-1}) \frac{1}{w_{ab}-1} \nonumber
        \label{eq:rbqn}
\end{eqnarray}

\subsection{In representation for $\fg=osp(1|2)$}

The universal $R$-matrix of $\cU_q(osp(1|2))$ was initially obtained in
\cite{KulRes,Sal}.  As indicated in Section \ref{sect:general}, we use the
finite $R$-matrix obtained as the evaluation of the universal $R$-matrix
given in \cite{KT141} (changing $q$ to $q^{-1}$):
\begin{equation*}
  \widetilde{R}_{12} = \left(
    \begin{array}{ccccccccc}
      -1 & 0 & 0\ &\ 0\ &\ 0\ &\ 0 & 0 & q^2-1 & 0 \\
      0 & 1 & 0 & 0 & 0 & 0 & 0 & 0 & 0 \\
      0 & 0 & 1 & 0 & 0 & 0 & q^{-1}-q \hspace{-1.4ex} & 0 & 0 \\
      0 & q^{-1}-q & 0 & 1 & 0 & 0 & 0 & 0 & 0 \\
      0 & 0 & 0 & 0 & q^{-1} & 0 & 0 & 0 & 0 \\
      q^{-1}-q & 0 & 0 & 0 & 0 & q & 0 &
      \hspace{-1.4ex} (q^{-1}-q)(q+1) \hspace{-1.4ex} & 0 \\
      0 & 0 & 0 & 0 & 0 & 0 & 1 & 0 & 0 \\
      0 & 0 & 0 & 0 & 0 & 0 & 0 & q & 0 \\
      0 & 0 & 0 & 0 & 0 & 0 & 0 & 0 & q^{-1} \\
    \end{array}
  \right)
\end{equation*}
The twist from $\uq{q}(osp(1|2))$ to ${\cal B}_{q,\lambda}(osp(1|2))$ is
represented by
\begin{eqnarray}
  F &=& \un \otimes \un - \frac{w(q-q^{-1})}{w-q} \; (E_{13} \otimes
  E_{31} - qE_{13} \otimes E_{12}) \nonumber \\
  && \;\; - \frac{qw(q-q^{-1})}{qw-1} \; (E_{21} \otimes E_{12} -
  q^{-1}E_{21} \otimes E_{31}) - \frac{w^2(q-q^{-1})(q+1)}{(qw-1)(w+1)}
  \; (E_{23} \otimes E_{32})
  \label{eq:twuqbqlosp12}
\end{eqnarray}
The resulting $R$-matrix of ${\cal B}_{q,\lambda}(osp(1|2))$ is given by
$F_{21}RF_{12}^{-1}$. The corresponding matrix $\widetilde{R}$ has
the following expression
\begin{equation*}
  \widetilde{R} =
  \left(
    \begin{array}{ccccccccc}
      r_{1111}& 0 & 0 & 0 & 0 & \frac{(q^2-1)qw}{q-w} & 0 &
      r_{1132} & 0
      \\
      0 & \frac{(q^3w-1)(w-q)}{q(qw-1)^2} & 0 &
      \frac{(q^2-1)w}{1-qw} & 0
      & 0 & 0 & 0 & 0 \\
      0 & 0 & 1 & 0 & 0 & 0 & \frac{q^2-1}{w-q} & 0 & 0 \\
      0 & \frac{q-q^{-1}}{qw-1} & 0 & 1 & 0 & 0 & 0 & 0 & 0 \\
      0 & 0 & 0 & 0 & q^{-1} & 0 & 0 & 0 & 0 \\
      \frac{q-q^{-1}}{qw-1} & 0 & 0 & 0 & 0 & q & 0 &
      \frac{(q^2-1)(q+1)}{(w-q)(w+1)} & 0 \\
      0 & 0 & \frac{(q-q^{-1})w}{q-w} & 0 & 0 & 0 &
      \frac{(q^3-w)(1-qw)}{q(q-w)^2} & 0 & 0 \\
      r_{3211} & 0 & 0 & 0 & 0 &
      \frac{(q^2-1)(q+1)w^2}{(1-qw)(w+1)} & 0
      & \frac{(q^2w+1)(q^2+w)}{q(w+1)^2} & 0 \\
      0 & 0 & 0 & 0 & 0 & 0 & 0 & 0 & q^{-1}
    \end{array}
  \right)
\end{equation*}
where $r_{1111} = \frac{q^2(w-1)^2+w(q-1)(q^3-1)}{q(qw-1)(q-w)}$, $r_{1132}
= \frac{(q^2-1)(q^2+w)(1-qw)}{(q-w)^2(w+1)}$ and $r_{3211} =
\frac{(q^2-1)w(q^2w+1)(w-q)}{q^2(qw-1)^2(w+1)}$.  It obeys the ordinary
dynamical Yang--Baxter equation.

\section{Twist from $\uq{}(\fg)$ to $\uq{s}(\fg)$\label{sect:Us}}
\setcounter{equation}{0}

\subsection{Universal form and cocycle condition}

We now construct the twist from $\uq{}(\fg)$ to $\uq{s}(\fg)$ as a limit
$q\rightarrow 1$, (i.e. $\hbar\rightarrow 0$ with $q=e^\hbar$) of the twist
(\ref{eq:FsommeBql}).  The consistency of this procedure follows from the
well-known Hopf algebra identification
$\uq{\hbar}(\fg)/(\hbar\;\uq{\hbar}(\fg)) \simeq \uq{}(\fg)$ \cite{Dri}.

In this quotient, the twist $\cF$ coincides with $\widehat{\cF}$.  It is
given by a formula analogous to (\ref{eq:FsommeBql}):
\begin{equation}
  \label{eq:FsommeUs}
  \cF_{[\uq{}(\fg) \to \uq{s}(\fg)]} = \un\otimes\un +
  \sum_{\{p,r\}\in (\cZ^*)^2} \varphi_{pr}(s) \; {\bold e}^p \otimes
  {\bold f}^r
\end{equation}
the functions $\varphi_{pr}(s)$ being the representatives of those
appearing in (\ref{eq:FsommeBql}).

It follows from the \emph{Hopf} algebra identification that
$\cF_{[\uq{}(\fg) \to \uq{s}(\fg)]}$ obeys in $\uq{}(\fg)^{\otimes 2}$ the
cocycle equation
\begin{equation}
  \label{eq:cocycleUs}
  \cF_{12}(s) (\Delta \otimes \un) \cF(s) = \cF_{23}(s+{h^\vee}^{(1)})
  (\un\otimes\Delta)  \cF(s)
\end{equation}
defining a \QTQHA\ denoted $\uq{s}(\fg)$ with $R$-matrix
$\cR(s) = \cF_{21}(s) \cF_{12}(s)^{-1}$.

The functions $\varphi_{pr}(s)$ satisfy recursion equations obtained as the
leading order in $\hbar$ of (\ref{eq:recurvarphiql}).  This procedure is
well-defined since the coefficient of $\varphi$ in the left hand side of
(\ref{eq:recurvarphiql}) is of order $\hbar$ and the coefficients in the
right hand side are at least of order 1 in $\hbar$ (due to the presence of
$\sigma_m$).  The leading order in $\hbar$ of (\ref{eq:recurvarphiql}) can
be expressed as the following equation in $\cF$:
\begin{equation}
  \label{eq:eqlinUs}
  [\cX,\cF] = \cF \; \widehat{r}
\end{equation}
where $\cX$ was defined in (\ref{eq:phiBql}) and
$\widehat{\cR}=\un\otimes\un + \hbar \, \widehat{r} + o(\hbar)$.  Note that
the equation (\ref{eq:eqlinUs}) is also the first non trivial term in the
expansion in $\hbar$ of (\ref{eq:eqlinBql}).

Under a similar hypothesis on $\cF$ as in Section \ref{sect:Bql}, replacing
$q$-deformed enveloping algebras by classical enveloping algebras, equation
(\ref{eq:eqlinUs}) has a unique solution expressed either by
(\ref{eq:FsommeUs}) or as the infinite product
\begin{equation}
  \label{eq:FprodUs}
  \cF_{[\uq{}(\fg) \to \uq{s}(\fg)]} = \prod_{k \ge
    1}^{\curvearrowleft} (\cX \otimes \un)^{-k} \;
  \Big(
    \un \otimes \un + (\cX \otimes \un)^{-1} \;\widehat{r}
  \Big)
  (\cX \otimes \un)^k
\end{equation}

\subsection{In representation for $\fg=sl_N$}

In the fundamental representation for $\fg=sl_N$, the evaluated
infinite product expression for $\cF$ reads
\begin{equation}
  F_{[\uq{}(sl_{N}) \to \uq{s}(sl_{N})]} = \prod_{k \ge
    1}^{\curvearrowleft} (X \otimes \un)^{-k} \; Y \; (X \otimes \un)^k
\end{equation}
with
\begin{eqnarray}
  && X = \sfrac{N-1}{N} \; \un + \diag(x_a)
  \qmbox{where}
  x_{a} = 2 s_{a} - 2 s_{a-1},\  s_{0} = s_{N} = 0 \;,
  \\
  && Y = \un \otimes \un + (X \otimes \un)^{-1}\; \widehat{r} =
  \un \otimes \un + (X
  \otimes \un)^{-1} \sum_{a<b} - 2 E_{ab} \otimes E_{ba} \;,
\end{eqnarray}
$\widehat{r} = - 2 \sum_{a<b} E_{ab} \otimes E_{ba}$ being the
classical $\widehat{R}$ matrix of $\uq{}(sl_{N})$.
\\
Computing the product leads to
\begin{equation}
  F = \un \otimes \un - \sum_{a<b} \frac{2}{x_{a}-x_{b}} \; E_{ab}
  \otimes E_{ba}
  \label{eq:freprsln}
\end{equation}
Regarding the $R$-matrix of $\uq{s}(sl_{N})$, its non
vanishing matrix elements
$R_{i_{1}i_{2}}^{j_{1}j_{2}}$ are given by ($1 \le a,b \le N$)
\begin{eqnarray}
        && R_{aa}^{aa} = 1 \nonumber \\
        && R_{ab}^{ab} = \left\{
        \begin{array}{ll}
                \displaystyle 1 & \qmbox{if}  b>a \\
                \displaystyle 1 - \frac{4}{(x_{a}-x_{b})^2}  &
                \qmbox{if}  b<a \\
        \end{array}
        \right.  \\
        && R_{ab}^{ba} = \frac{2}{x_{a}-x_{b}} \qmbox{for} a\neq b \;
        \nonumber
        \label{eq:rUs}
\end{eqnarray}
which indeed satisfies the dynamical Yang--Baxter equation (\ref{eq:DYBE0}).

\subsection{In representation for $\fg=osp(1|2)$}

Using again the infinite product expression (\ref{eq:FprodUs}),
the twist from $\uq{}(osp(1|2))$ to $\uq{s}(osp(1|2))$ is represented by
\begin{equation}
  F = \un \otimes \un - \frac{2}{s-1} \; (E_{13} \otimes E_{31} - E_{13}
  \otimes E_{12}) - \frac{2}{s+1} \; (E_{21} \otimes E_{12} + E_{23}
  \otimes E_{32} - E_{21} \otimes E_{31})
  \label{eq:twuusosp12}
\end{equation}
The $R$-matrix of $\uq{s}(osp(1|2))$ is then straightforwardly $R(s) =
F_{21}(s) F_{12}^{-1}(s)$.  The corresponding matrix $\widetilde{R}$
satisfies the ordinary dynamical Yang--Baxter equation.

\subsection{In representation for $\fg=sl(1|2)$}

Similarly,
the twist from $\uq{}(sl(1|2))$ to $\uq{s}(sl(1|2))$ is represented by
\begin{equation}
  F = \un \otimes \un + \frac{1}{s_2-1} \; E_{21} \otimes E_{12} 
  - \frac{1}{s_1+s_2} \; E_{31} \otimes E_{13} 
  + \frac{1}{s_1+1} \; E_{32} \otimes E_{23}
  \label{eq:twuussl12}
\end{equation}
The $R$-matrix of $\uq{s}(sl(1|2))$ is then given by
$R(s) = F_{21}(s) F_{12}^{-1}(s)$.  
Again, the corresponding matrix $\widetilde{R}$
satisfies the ordinary dynamical Yang--Baxter equation.

\section{Twist from $\uq{q}(\fgh)$ to $\uq{q,\lambda}(\fgh)$}
\setcounter{equation}{0}

\subsection{Universal form}

The expression of a universal twist from $\uq{q}(\fg)$ to
$\cB_{q,\lambda}(\fg)$ was recalled in Section \ref{sect:Bql}.  Using the
fact that $\uq{q}(\fg)$ is a Hopf subalgebra of $\uq{q}(\fgh)$, we use this
twist to construct the dynamical $\uq{q,\lambda}(\fgh)$.  Indeed,
$\cF_{[\uq{q}(\fg) \to {\cal B}_{q,\lambda}(\fg)]}$ seen as an element of
$\uq{q}(\fgh)^{\otimes 2}$ satisfies the shifted cocycle condition,
yielding a dynamical $R$-matrix
\begin{equation}
  \label{eq:twistUql}
  \cR_{\uq{q,\lambda}(\fgh)}(w) =
  \cF_{21}(w) \; \cR_{\uq{q}(\fgh)} \; \cF_{12}^{-1}(w) \;.
\end{equation}
This defines $\uq{q,\lambda}(\fgh)$ as a \QTQHA. We now identify it, by
computing its evaluation representation for $\fg=sl_N$, to a realisation of
the $p \to 0$ limit of $\elpb(\sln)$ defined in Appendix B.

\subsection{In representation for $\fg=sl_N$}

One gets the matrix elements $R_{i_{1}i_{2}}^{j_{1}j_{2}}$ of the
$R$-matrix of $\uq{q,\lambda}(\sln)$ ($1 \le a,b \le N$):
\begin{eqnarray}
  && R_{aa}^{aa} = \rho_{\uq{q,\lambda}}(z)  \nonumber \\
  && R_{ab}^{ab} = \rho_{\uq{q,\lambda}}(z) \left\{
    \begin{array}{ll}
      \displaystyle \frac{q(1-z)}{1-q^2z} & \mbox{if} \;\; b>a \\
      & \\
      \displaystyle

\frac{q(1-z)}{(1-q^2z)}\,\frac{(1-w_{ab}q^2)(1-w_{ab}q^{-2})}{(1-w_{ab})^2}
      & \mbox{if} \;\; b<a \\
    \end{array}
  \right.  \nonumber\\[3mm]
  && R_{ab}^{ba} = \rho_{\uq{q,\lambda}}(z) \;
  \frac{(1-q^2)(1-w_{ab}z)}{(1-q^2z)(1-w_{ab})}
  \label{eq:ruqln}
\end{eqnarray}
the normalisation factor being given by
\begin{equation}
  \rho_{\uq{q,\lambda}}(z) = q^{-\frac{N-1}{N}} \;
    \frac{(q^{2}z;q^{2N})_{\infty} \;
    (q^{2N-2}z;q^{2N})_{\infty}}{(z;q^{2N})_{\infty} \;
    (q^{2N}z;q^{2N})_{\infty}} \;.
  \label{eq:rhoz}
\end{equation}
This $R$-matrix satisfies the Dynamical Yang--Baxter equation
(\ref{eq:DYBE}).  It is indeed the limit $p \to 0$ of the $R$-matrix
(\ref{eq:rbqpnbis}).

\subsection{In representation for $\fg=osp(1|2)$}

We first construct a represented $R$-matrix of $\widehat{osp(1|2)}$ through
a Baxterisation procedure \cite{Jones}.  We get two $R$-matrices with
spectral parameter constructed from $\widetilde{R}_{12}$,
$\widetilde{R}_{21}^{-1}$ and $\widetilde{P}_{12}$ (the non-graded
permutation) which obey the non-graded Yang--Baxter equation (with spectral
parameter).  They read:
\begin{equation}
\widetilde{R}(z)=\frac{1-z}{(1-za)(1-zq^2)}\,\widetilde{R}_{{12}} -
\frac{aq^2z(1-z)}{(1-za)(1-zq^2)}\,\widetilde{R}_{21}^{-1} +
\frac{z(1-a)(1-q^2)}{(1-za)(1-zq^2)}\,\widetilde{P}_{12} \;\;
\mbox{with} \;\; a=-q
\;\; \mbox{or} \;\; q^3
\end{equation}
where the normalisations are such that $\widetilde{R}(0) =
\widetilde{R}_{12}$, $\widetilde{R}(\infty) =
q^{-2}\widetilde{R}_{21}^{-1}$ and $\widetilde{R}(1) =
q^{-1}\widetilde{P}_{12}$.
Explicitly, for $a=-q$, we
get:
\begin{equation*}
  \widetilde{R}(z)=\left(
        \begin{array}{ccccccccc}
          r_{1111}& 0 & 0 & 0 & 0 & \frac{z(1-z)(q^2-1)}{(1-zq^2)(1+zq)} & 0 &
          \frac{(1-z)(q^2-1)}{(1-zq^2)(1+zq)} & 0\\
          0 &\frac{1-z}{1-zq^2} & 0 & \frac{(q^{-1}-q)z}{1-zq^2} & 0 & 0 & 0 &
          0 & 0 \\
          0 & 0 & \frac{1-z}{1-zq^2} & 0 & 0 & 0 & \frac{q^{-1}-q}{1-zq^2} & 0
          & 0 \\
          0 & \frac{q^{-1}-q}{1-zq^2} & 0 & \frac{1-z}{1-zq^2} & 0 & 0 & 0 & 0
          & 0 \\
          0 & 0 & 0 & 0 & q^{-1} & 0 & 0 & 0 & 0 \\
          \frac{(1-z)(q^{-1}-q)}{(1-zq^2)(1+zq)} & 0 & 0 & 0 & 0 &
          \frac{(1-z)(z+q)}{(1-zq^2)(1+zq)} & 0 &
          \frac{(q^{-1}-q)(q+1)}{(1-zq^2)(1+zq)} & 0 \\
          0 & 0 & \frac{(q^{-1}-q)z}{1-zq^2} & 0 & 0 & 0 &
\frac{1-z}{1-zq^2} &
          0 & 0 \\
          \frac{z(1-z)(q^{-1}-q)}{(1-zq^2)(1+zq)} & 0 & 0 & 0 & 0 &
          \frac{z^2(q^{-1}-q)(q+1)}{(1-zq^2)(1+zq)} & 0 &
          \frac{(1-z)(z+q)}{(1-zq^2)(1+zq)} & 0 \\
          0 & 0 & 0 & 0 & 0 & 0 & 0 & 0 & q^{-1}
        \end{array}
  \right)
\end{equation*}
with $r_{1111}=\frac{z-1}{1-zq^2}+\frac{z(q+1)(q^{-1}-q)}{(1-zq^2)(1+zq)}$.
\\
Applying to the above matrix the twist given by formula
(\ref{eq:twuqbqlosp12}) allows one to get the explicit evaluated $R$-matrix
for the \QTQHA\ $\uq{q,\lambda}(osp(1|2))$.  Due to its cumbersome nature,
we omit this explicit form here.

\section{Twist from $\dy{}{}(\fg)$ to $\dy{s}{}(\fg)$}
\setcounter{equation}{0}

\subsection{Universal form}

We similarly use the fact that $\uq{}(\fg )$ is a Hopf subalgebra of
$\dy{}{}(\fg )$, $\fg$ belonging to the (super) unitary series.  Hence
the cocycle identity (\ref{eq:cocycleUs}) is also 
an identity in $\dy{}{}(\fg )^{\otimes 3}$ for $\cF$ defined as in
(\ref{eq:FsommeUs}), now considered as an element of $\dy{}{}(\fg
)^{\otimes 2}$.  This twist, applied to the universal $R$-matrix of
$\dy{}{}(\fg)$ (given in \cite{KT}), yields a dynamical $R$-matrix
\begin{equation}
  \label{eq:Rdys}
  \cR(s)=\cF_{21}(s)\, \cR\, \cF_{12}^{-1}(s)
\end{equation}
which characterise $\dy{s}{}(\fg )$ as a \QTQHA.

\subsection{In representation for $\fg=sl_N$}

We apply the twist (\ref{eq:freprsln}) to the $R$-matrix of
$\dy{}{}(sl_{N})$, given in (\ref{eq:rdy})-(\ref{eq:rhody}), to get the
$R$-matrix of $\dy{s}{}(sl_{N})$.  In the fundamental representation, it
has the following non vanishing elements $R_{i_{1}i_{2}}^{j_{1}j_{2}}$ ($1
\le a,b \le N$)
\begin{eqnarray}
  && R_{aa}^{aa} = \rho_{DY_{s}}(u)  \nonumber \\
  && R_{ab}^{ab} = \rho_{DY_{s}}(u) \left\{
    \begin{array}{ll}
      \displaystyle \frac{u}{u+1} & \mbox{if} \;\; b>a \\
      & \\
      \displaystyle \left( 1 - \frac{4}{(x_{a}-x_{b})^2} \right) \;
      \frac{u}{u+1} & \mbox{if} \;\; b<a \\
    \end{array}
  \right.  \nonumber \\
  && R_{ab}^{ba} = \rho_{DY_{s}}(u) \left( 1 + \frac{2u}{x_{a}-x_{b}}
  \right) \;
  \frac{1}{u+1}
  \label{eq:rdys}
\end{eqnarray}
the normalisation factor being $\rho_{DY_{s}}(u) = \rho_{\dy{}{}}(u)$.
This matrix can also be obtained as the scaling limit of the
$R$-matrix (\ref{eq:ruqln}) of $\cU_{q,\lambda}(\sln)$.
It satisfies the dynamical Yang--Baxter equation (\ref{eq:DYBE2}).

\subsection{In representation for $\fg=sl(1|2)$}

Similarly, by applying the twist (\ref{eq:twuussl12}) to the $R$-matrix
of $\dy{}{}(sl(1|2))$, one gets the evaluated $R$-matrix for the \QTQHA\
$\dy{s}{}(sl(1|2))$:
\begin{equation*}
  \widetilde{R} = \left(
    \begin{array}{ccccccccc}
      \frac{1-u}{1+u} & 0 & 0 & 0 & 0 &
      0 & 0 & 0
      & 0 \\
      0 & \frac{u s_2 (s_2-2)}{(1+u)(s_2-1)^2} & 0 &
      \frac{s_2-1+u}{(1+u)(s_2-1)} & 0 & 0 & 0 & 0 & 0 \\
      0 & 0 & \frac{u(1-(s_1+s_2)^2)}{(1+u)(s_1+s_2)^2} & 0 & 0 & 0 &
      \frac{s_1+s_2+u}{(1+u)(s_1+s_2)} & 0 &
      0 \\
      0 & \frac{s_2-1-u}{(1+u)(s_2-1)} & 0 & \frac{u}{1+u} & 0 & 0 &
      0 & 0 &
      0 \\
      0 & 0 & 0 & 0 & 1 & 0 & 0 & 0 & 0 \\
      0 & 0 & 0 & 0 & 0 & \frac{u s_1 (s_1+2)}{(1+u)(s_1+1)^2} & 0 &
      \frac{s_1+1+u}{(1+u)(s_1+1)} & 0 \\
      0 & 0 & \frac{s_1+s_2-u}{(1+u)(s_1+s_2)} & 0 & 0 & 0 &
      \frac{-u}{1+u} & 0 & 0 \\
      0 & 0 & 0 & 0 & 0 &
      \frac{s_1+1-u}{(1+u)(s_1+1)} & 0 & \frac{u}{1+u} & 0 \\
      0 & 0 & 0 & 0 & 0 & 0 & 0 & 0 & \frac{1-u}{1+u}
    \end{array}
  \right)
\end{equation*}
It satisfies the dynamical Yang--Baxter equation (\ref{eq:DYBE2}).

\appendix
\section{Notations}
\setcounter{equation}{0}

Multiple Gamma functions are defined by
\begin{equation}
  \Gamma_{r}(x \vert \omega_{1},\ldots,\omega_{r}) = \exp\left(
  \frac{\partial}{\partial s} \sum_{n_{1},\ldots,n_{r} \ge 0}
  (x+n_{1}\omega_{1}+\ldots+n_{r}\omega_{2})^{-s}\Bigg\vert_{s=0} \right) \,.
\end{equation}
Barnes' multiple sine function  $S_{r}(x \vert
\omega_{1},\ldots,\omega_{r})$ of periods $\omega_{1}$, $\ldots$,
$\omega_{r}$ is defined by \cite{Barnes}
\begin{equation}
  S_{r}(x \vert \omega_{1},\omega_{2}) =
  \Gamma_{r}(\omega_{1}+\ldots+\omega_{r}-x \vert
  \omega_{1},\ldots,\omega_{r})^{(-1)^r} \Gamma_{r}(x \vert
  \omega_{1},\ldots,\omega_{r})^{-1} \;,
\end{equation}
They satisfy for each $i \in [1,\ldots,r]$
\begin{eqnarray}
  &&
  \frac{\Gamma_{r}(x+\omega_{i} \vert
    \omega_{1},\ldots,\omega_{r})}
  {\Gamma_{r}(x
    \vert \omega_{1},\ldots,\omega_{r})} = \frac{1}{\Gamma_{r-1}(x \vert
    \omega_{1},\ldots,\omega_{i-1},\omega_{i+1},\ldots,\omega_{r})} \,,
  \\ &&
  \frac{S_{r}(x+\omega_{i} \vert \omega_{1},\ldots,\omega_{r})}{S_{r}(x
    \vert \omega_{1},\ldots,\omega_{r})} = \frac{1}{S_{r-1}(x \vert
    \omega_{1},\ldots,\omega_{i-1},\omega_{i+1},\ldots,\omega_{r})} \,.
\end{eqnarray}
In particular,
\begin{equation}
  \Gamma_1(x|\omega_1) = \frac{\omega_1^{x/\omega_1}}{\sqrt{2\pi \omega_1}}
  \; \Gamma\left(\frac{x}{\omega_1}\right)\;,
  \qquad\qquad
  S_1(x|\omega_1) = {2\sin\frac{\pi x}{\omega_{1}}} \,.
\end{equation}

\section{Definition of $\elpb(\sln)$}
\setcounter{equation}{0}

The quantum affine elliptic algebra $\elpb(\sln)$ was originally defined
in \cite{Fe} using the $RLL$ formalism.  The characteristic $R$-matrix
takes the following form ($1 \le a,b \le N$):
\begin{eqnarray}
        R(z) &=& \rho_{\elpb}(z) \left( \sum_{a} E_{aa} \otimes E_{aa} +
\sum_{a
        \ne b} q^2 \frac{\Theta_{p}(q^{-2}w_{ab})}{\Theta_{p}(w_{ab})} \;
        \frac{\Theta_{p}(z)}{\Theta_{p}(q^2z)} \; E_{aa} \otimes E_{bb}
\right.
        \nonumber \\
        && + \; \left.  \sum_{a \ne b}
        \frac{\Theta_{p}(w_{ab}z)}{\Theta_{p}(w_{ab})} \;
        \frac{\Theta_{p}(q^2)}{\Theta_{p}(q^2z)} \; E_{ab} \otimes E_{ba}
        \right)
        \label{eq:rbqpn}
\end{eqnarray}
where $\Theta_p(z) = (z;p)_\infty \, (pz^{-1};p)_\infty \, (p;p)_\infty$
and $w_{ab} = q^{x_{a}-x_{b}}$, the infinite multiple products being
defined by $(z;p_1,\dots,p_m)_\infty = \prod_{n_i \ge 0} (1-z p_1^{n_1}
\dots p_m^{n_m})$.  The normalisation factor $\rho_{\elpb}(z)$ is given by
\begin{equation}
  \rho_{\elpb}(z) = q^{-\frac{N-1}{N}} \; \frac{(q^2z;q^{2N},p)_{\infty} \;
  (q^{2N-2}z;q^{2N},p)_{\infty}} {(z;q^{2N},p)_{\infty} \;
  (q^{2N}z;q^{2N},p)_{\infty}} \; \frac{(pz^{-1};q^{2N},p)_{\infty} \;
  (pq^{2N}z^{-1};q^{2N},p)_{\infty}} {(pq^2z^{-1};q^{2N},p)_{\infty} \;
  (pq^{2N-2}z^{-1};q^{2N},p)_{\infty}} \;.
  \label{eq:rhoelpb}
\end{equation}
It was proven in \cite{JKOS} that this quantum affine elliptic algebra was
a \QTQHA\ obtained by a Drinfel'd twist of shifted-cocycle type from
$\uq{q}(\sln)$.  The $R$-matrix thus obtained is actually a gauge transform
of (\ref{eq:rbqpn}) where the coefficients of $E_{aa} \otimes E_{bb}$
become
\begin{eqnarray}
        && R_{ab}^{ab} = \rho_{\elpb}(z)\left\{
        \begin{array}{ll}
                \displaystyle q \; \frac{(pw_{ab}^{-1}q^2;p)_{\infty} \;
                (pw_{ab}^{-1}q^{-2};p)_{\infty}}
{(pw_{ab}^{-1};p)_{\infty}^2} \;
                \frac{\Theta_{p}(z)}{\Theta_{p}(q^2z)} & \mbox{if} \;\; b>a \\
                & \\
                \displaystyle q \; \frac{(w_{ab}^{-1}q^2;p)_{\infty} \;
                (w_{ab}^{-1}q^{-2};p)_{\infty}}
{(w_{ab}^{-1};p)_{\infty}^2} \;
                \frac{\Theta_{p}(z)}{\Theta_{p}(q^2z)} & \mbox{if} \;\; b<a \\
        \end{array}
        \right.
        \label{eq:rbqpnbis}
\end{eqnarray}

\bigskip

\paragraph{Acknowledgments:} 

This work was supported in part by CNRS
and EC network contract number FMRX-CT96-0012. 
J.A.  wishes to thank the LAPTH for its kind hospitality. 
The authors thank M.~Rossi for discussions at early stage of this
work.


\end{document}